\documentclass{IEEEtran}

\usepackage{graphicx}
\usepackage{amsmath}
\usepackage{amssymb}
\usepackage{amsfonts}
\usepackage{mathrsfs}
\usepackage{dsfont}
\makeatletter
\def\amsbb{\use@mathgroup \M@U \symAMSb}
\makeatother
\usepackage{subfigure}
\usepackage[usenames,dvipsnames]{xcolor}
\usepackage{epstopdf}

\allowdisplaybreaks
\usepackage{balance}

\newtheorem{theorem}{Theorem}
\newtheorem{remark}{Remark}
\newtheorem{definition}{Definition}
\newtheorem{lemma}{Lemma}
\newtheorem{proposition}{Proposition}

\newtheorem{assumption}{Assumption}
\usepackage{float,graphicx,subfig,tikz}
\usetikzlibrary{shapes,arrows,scopes,plotmarks}

\newcommand{\vect}[1]{\boldsymbol{#1}} 

\usepackage{enumerate}
\usepackage[prependcaption,colorinlistoftodos]{todonotes}

\newcommand{\n}[1]{\textcolor{black}{#1}}
\usepackage[normalem]{ulem}

\title{\color{black}Stability and optimality of distributed secondary frequency control schemes in power networks}

\author{Andreas Kasis\thanks{Andreas Kasis, {\color{black}Nima Monshizadeh} and Ioannis Lestas are with the Department of Engineering, University of Cambridge, Trumpington Street, Cambridge, CB2 1PZ, United Kingdom; e-mails: ak647@cam.ac.uk, {\color{black}n.monshizadeh@eng.cam.ac.uk}, icl20@cam.ac.uk}, {\color{black}Nima Monshizadeh}, Eoin Devane\thanks{Eoin Devane is with the Cambridge Centre for Analysis, Centre for Mathematical Sciences, University of Cambridge, Wilberforce Road, Cambridge, CB3 0WA, United Kingdom; e-mail: esmd2@cam.ac.uk},  and Ioannis Lestas
\thanks{\color{black}Paper \cite{KDI_secondary_CDC} is a conference version of this manuscript. {\color{black}This manuscript includes proofs of our main results, additional discussion and {\color{black}additional} results of independent interest.}}}

\begin{document}
\tikzstyle{block} = [draw, fill=white!20, rectangle,
    minimum height=3em, minimum width=6em]
\tikzstyle{sum} = [draw, fill=blue!20, circle, node distance=1cm]
\tikzstyle{input} = [coordinate]
\tikzstyle{output} = [coordinate]
\tikzstyle{pinstyle} = [pin edge={to-,thin,black}]

\tikzstyle{int}=[draw, fill=blue!20, minimum size=2em]
\tikzstyle{init} = [pin edge={to-,thin,black}]

\maketitle

\begin{abstract}
We present a \n{systematic} method {\color{black} for designing} distributed generation and demand control schemes for secondary frequency regulation in power networks {\color{black}such that stability and {an economically optimal power} allocation can be guaranteed}. A dissipativity condition is imposed on net power supply variables {\color{black} to provide stability guarantees. Furthermore, economic optimality is achieved by explicit {\color{black}decentralized} steady state conditions on {\color{black}the} generation and controllable {\color{black}demand.}}
 We discuss how various classes of dynamics used in recent studies fit within our framework and give examples of higher order generation and controllable demand dynamics that can be included within our analysis. \n{In case of linear dynamics, {\color{black}we discuss how} the proposed dissipativity condition can be {\color{black}efficiently} verified using an appropriate linear matrix inequality.}
{\color{black}Moreover, it  is shown how the addition of a suitable observer layer can relax the requirement for demand measurements in the employed controller.}
\n{The efficiency and practicality of the proposed results are demonstrated {\color{black}with a simulation} on the Northeast Power Coordinating Council (NPCC) 140-bus {\color{black}system.}}
\end{abstract}

\vspace{-2mm}
\section{Introduction}

{\color{black}
Renewable sources of energy are expected to grow in penetration within power networks over the next years \cite{Lund, Future_Grid}.
Moreover, {\color{black}it is anticipated that} controllable loads will  be incorporated within power networks in order to provide benefits such as fast response to changes in power generated from renewable sources and the ability for peak demand reduction.
Such changes will greatly increase power network complexity revealing a need for highly distributed schemes that will guarantee its stability when `plug and play' devices are incorporated.
In the recent years, research attention has {\color{black}increasingly} focused on such distributed schemes with studies regarding both primary (droop) control as in ~\cite{Primary2, trip_depersis, kasis_dev_lestas_1} and secondary control as in \cite{P30_Low}, \cite{P35_DePersis}.
}


An issue of economic optimality in {\color{black}the} power allocation is raised if
{highly distributed schemes are to be used for frequency control.}
Recent studies attempted to address this issue \n{by crafting the equilibrium of the system such that it coincides with
the optimal solution of a suitable network optimization problem.}
{\color{black}To establish optimality of an equilibrium in a distributed fasion, it is evident that a synchronising variable is required.} While in the primary control, frequency is used as the synchronising variable (e.g. \cite{kasis_dev_lestas_1,P19,P34 - Low,kasis_dev_lestas_2}), in the secondary control
{\color{black} a different variable is {\color{black}synchronized} by making use of  information exchanged between buses}
\cite{P30_Low,P35_DePersis,LC3,P21}.

{\color{black}Over the last few years many studies have attempted to address issues regarding stability and optimization in secondary frequency control.}
An important feature in many of those is that the dynamics {\color{black} considered follow from a primal/dual
algorithm associated with some optimal power allocation problem} \cite{P30_Low}, \cite{P6_Low}, \cite{P57_Low}, \cite{P1_Low}.
This is a powerful approach that reveals the information structure needed to achieve optimality and satisfy the {\color{black}constraints involved}.
{\color{black}Nevertheless, when higher order generation dynamics need to be considered, these do not necessarily follow as gradient dynamics of a corresponding optimization problem and therefore alternative approaches need to be employed.}

{\n{Another trend in the secondary frequency control is the use of distributed averaging proportional integral (DAPI) controllers
{\cite{Dorfler_DAPI, Bullo_DAPI, Dorfler_DAPI_2,DePersis_automatica,Martin_ACC}}. Advantages of DAPI controllers lie in their simplicity as they only measure local frequency and {\color{black}exchange} a synchronization signal in a distributed fashion {\color{black}without requiring load and power flow measurements}. On the other hand, it is not easy to accommodate line and power flow constraints, and higher-order generation and controllable demand dynamics in this setting. Moreover the existing results in this context are limited to the case of proportional active power sharing and quadratic cost functions.}

{\color{black}One of our aims in this paper is to present a methodology that allows to incorporate}
general {\color{black}classes of higher order} generation and demand control dynamics {\color{black}while} \n{ensuring}  stability and optimality of the equilibrium points.
Our analysis borrows ideas from our previous work in \cite{kasis_dev_lestas_1} 
and adapts those to secondary frequency control, by incorporating the additional communication layer needed in this context. In particular, we consider general classes of aggregate power supply dynamics at each bus and impose two conditions; a dissipativity condition that ensures stability, and a \n{steady-state} condition that ensures optimality of the power allocation.  \n{An important feature of these conditions {\color{black}is} that they {\color{black}are} decentralized.}
\n{{\color{black}Furthermore, in} the case of linear supply dynamics, the proposed dissipativity condition can be efficiently verified by means of a linear matrix inequality (LMI).}
Various examples are also described to illustrate the significance of our approach and the way it could facilitate \n{a systematic analysis and} design. {
{\color{black} Finally, we discuss how an appropriately designed observer, allows to relax the requirement of an explicit knowledge of the uncontrollable demand, and show that the stability and optimality guarantees remain valid in this case.}

{\color{black}The paper is structured} as follows. {\color{black}Section \ref{sec: Preliminaries} provides some basic notation and preliminaries}. In section \ref{Problem_Formulation} we present the power network {\color{black}model}, the classes of generation and controllable demand dynamics and the {\color{black}optimization} problem to be considered. {\color{black}Sections \ref{Assumptions} and \ref{Main_Results} include our main assumptions and results.} In Section \ref{Discussion} we discuss how the results apply to various dynamics for generation and {\color{black}demand, provide intuition regarding our analysis and} show how the controller requirements may be relaxed by incorporating an appropriate observer. In section \ref{Simulation}, we demonstrate our results through a simulation on {\color{black}the NPCC 140-bus} system. Finally, conclusions are drawn in section \ref{Conclusion}.

\section{\color{black}Notation and Preliminaries}\label{sec: Preliminaries}

 Real numbers are denoted by $\mathbb{R}$, and the set of $n$-dimensional vectors with real entries is denoted by $\mathbb{R}^n$.
For a function $f(q)$, {\color{black}$f:\mathbb{R}\rightarrow \mathbb{R}$,} we denote its first derivative by $f'(q) = \tfrac{d}{dq} f(q)${\color{black}, its inverse by $f^{-1}(.)$. 
{\color{black}A function $f:\mathbb{R}^n \rightarrow \mathbb{R}$
is said to be positive
semidefinite
if $f(x)\geq0$.
It is positive definite
if $f(0) = 0$ and $f(x) > 0$ for every $x\neq0$.
We say that $f$ is positive definite with respect to component $x_j$ if {\color{black}$f(x) = 0$ {\color{black}implies} $x_j=0$,} and $f(x) > 0$ for every
$x_j \neq0$.}
{\color{black} A function $f:X\rightarrow Y$ is called surjective if $\forall y \in Y, \exists x \in X$ such that $f(x) = y$.}
 For $a,b \in \mathbb{R}$, $a\leq b$, the expression $[q]^b_a$ will be used to denote $\max\{\min\{q,b\},a\}$ and we write $\vect{0}_n$ to denote $n \times 1$ vector with all elements equal to $0$.
We use $\mathds{1}_{a \leq b}$ to denote a function that takes the value of $1$ when $a \leq b$, for $a,b \in \mathbb{R}$, and of $0$ otherwise. {\color{black}The Laplace transform of a signal $h(t)$, $h:\mathbb{R}\rightarrow \mathbb{R}$, is {denoted by $\hat{h}(s) = \int_0^\infty e^{-st} h(t) \, dt$.}}
{\color{black}Finally, for input/output systems $B_j, j=1,\hdots,N$, with respective inputs $u_j$ and outputs $y_j$,
their direct sum, denoted by $\bigoplus_{j=1}^N B_j$, represents a system with input $[u_1^T, u_2^T, \dots u_N^T]^T$ and output $[y_1^T, y_2^T, \dots y_N^T]^T$.}


{\color{black}Within the paper, we will consider {\color{black}subsystems\footnote{\color{black}Note that such subsystems will be used to characterize generation and demand dynamics and will be explicitly stated when considered.} that will be modeled as dynamical systems}
with input {\color{black}$u(t) \in \mathbb{R}^m$, state $x(t) \in \mathbb{R}^n$,}
and output $y(t) \in \mathbb{R}^k$ and a state space realization 
\vspace{-.1cm}
\begin{equation} \label{dynsys}
\begin{aligned}
&\dot{x} = f(x,u),\\
&y = g(x,u),
\end{aligned}
\end{equation}
where {\color{black}$f : \mathbb{R}^n \times \mathbb{R}^m \to \mathbb{R}^n$ is locally Lipschitz and $g : \mathbb{R}^n \times \mathbb{R}^m \to \mathbb{R}^k$ is continuous.
{\color{black}We assume in~\eqref{dynsys}
that given any constant input $u(t) \equiv \bar{u}$, there exists a unique\footnote{\color{black}The uniqueness assumption on the equilibrium point for a given input could be relaxed to having isolated equilibrium points, but it is used here for simplicity in the presentation.}
locally asymptotically stable equilibrium point $\bar{x} \in \mathbb{R}^m$, i.e. $f(\bar{x}, \bar{u}) = 0$.
The region of attraction\footnote{\color{black}That is, for the constant input $\zeta_j = \bar{\zeta}_j$, any solution $x(t)$ {\color{black}of~\eqref{sys2}} with initial condition $x(0) \in X_0$ must satisfy $x(t) \to \bar{x}$ as $t \to \infty$.}
of~$\bar{x}$ is denoted by $X_0$.}
We also define} the static input-state characteristic {\color{black} map $k_x : \mathbb{R}^m \to \mathbb{R}^n$ as
\begin{equation*} 
k_x(\bar{u}) := \bar{x},
\end{equation*}
and the static input-output characteristic map
$k_y : \mathbb{R}^m \to \mathbb{R}^k$,} 
\begin{equation}
k_y(\bar{u}) := g(k_x(\bar{u}), \bar{u}).
\label{iochar}
\end{equation}} }

\vspace{-2mm}
\section{Problem formulation} \label{Problem_Formulation}

\subsection{Network model}\label{Network_Model}

We describe the power network model by a connected graph $(N,E)$ where $N = \{1,2,\dots,|N|\}$ is the set of buses and $E \subseteq N \times N$ the set of transmission lines connecting the buses.
{\color{black}There are two types of buses in the network, buses with inertia and buses without inertia.
Since generators have inertia, it is reasonable to assume that only buses with inertia have non-trivial generation dynamics. We define $G = \{1,2,\dots,|G|\}$ and $L= \{|G|+1,\dots,|N|\}$ as the sets buses with and without inertia respectively such that $|G| + |L| = |N|$.}
Moreover, the term $(i,j)$ denotes the link connecting buses $i$ and $j$. The graph $(N,E)$ is assumed to be directed with \n{an} arbitrary direction, so that if $(i,j) \in E$ then $(j,i) \notin E$. Additionally, for each $j \in N$, we use $i:i\rightarrow j$ and $k:j\rightarrow k$ to denote the sets of buses that precede and succeed bus~$j$ respectively. It should be noted that the form of the dynamics in~\eqref{sys1}--\eqref{sys2} below is not affected by changes in graph ordering, and our results are independent of the choice of direction.  We make the following assumptions for the network: \newline
1) Bus voltage magnitudes are $|V_j| = 1$ p.u. for all $j \in N$.
2) Lines $(i,j) \in E$ are lossless and characterized by their susceptances $B_{ij} = B_{ji} > 0$. \newline
3) Reactive power flows do not affect bus voltage phase angles and frequencies. \newline
{\color{black}Such} assumptions are generally valid at medium to high voltages or when tight voltage control is present, and are often used in secondary frequency control studies
\cite{Bergen_Vittal}.

Swing equations can then be used to describe the rate of change of frequency at generation buses. Power must also be conserved at each of the load buses.
 This motivates the following system dynamics (e.g.~\cite{Bergen_Vittal}),
\begin{subequations} \label{sys1}
\begin{equation}
\dot{\eta}_{ij} = \omega_i - \omega_j, \; (i,j) \in E, \label{sys1a}
\end{equation}
\begin{equation}
M_j \dot{\omega}_j = - p_j ^L + p_j^M - (d^c_j + d^u_j) - \sum_{k:j\rightarrow k} p_{jk} + \sum_{i:i\rightarrow j} p_{ij}, \; j\in G, \label{sys1b}
\end{equation}
\begin{equation}
 0 = - p_j ^L - (d^c_j + d^u_j) - \sum_{k:j\rightarrow k} p_{jk} + \sum_{i:i\rightarrow j} p_{ij}, \; j\in L, \label{sys1c}
\end{equation}
\begin{equation}
 p_{ij}=B_{ij} \sin\eta _{ij} - p_{ij}^{nom}, \; (i,j) \in E. \label{sys1d}
\end{equation}
\end{subequations}
In system~\eqref{sys1}, the time-dependent variables $\omega_j$, $d^c_j$ and $p^M_j$ represent, respectively, deviations from a nominal value\footnote{A nominal value of a variable is defined as its value at an equilibrium of \eqref{sys1} with frequency at its nominal value of 50Hz (or 60Hz).} for the frequency and controllable load at bus $j$ and the mechanical power injection to the {\color{black} generation} bus $j$.
 The quantity $d^u_j$ represents the {\n{uncontrollable frequency-dependent}  load and generation  damping present at bus $j$.}
 The time-dependent variables $\eta_{ij}$ and $p_{ij}$ represent, respectively, the power angle difference\footnote{\color{black}The quantities $\eta_{ij}$ represent the phase differences between buses $i$ and $j$, given by $\theta_i - \theta_j$, i.e. $\eta_{ij} = \theta_i - \theta_j$. The angles themselves must also satisfy $\dot{\theta}_j = \omega_j$ at all $j \in N$. This equation is omitted in \eqref{sys1} since the power transfers are functions of the phase differences only.}
 and the deviation of the power transferred from bus $i$ to bus $j$ from the nominal value, $p_{ij}^{nom}$. The constant $M_j > 0$ denotes the generator inertia.  The response of the system~\eqref{sys1} will be studied, when a step change $p_j^L, j \in N$ {\color{black} occurs in} the uncontrollable demand.

\n{{\color{black}In order to investigate broad classes of generation and demand dynamics and control policies,} we let the scalar variables $p^M_j$, $d^c_j$, and $d^u_j$ be generated by dynamical systems of form \eqref{dynsys}, namely}
\begin{subequations} \label{sys2}
\begin{equation} \label{sys2p}
\begin{aligned}
&\dot{x}^{M,j} = f^{M,j}(x^{M,j},\zeta_j), \\
&p^M_j = g^{M,j}(x^{M,j},\zeta_j),
\end{aligned} \hspace{2em}j \in G,
\end{equation}
\begin{equation} \label{sys2dc}
\begin{aligned}
&\dot{x}^{c,j} = f^{c,j}(x^{c,j},\zeta_j), \\
&d^c_j = g^{c,j}(x^{c,j},\zeta_j),
\end{aligned} \hspace{2em}j \in N,
\end{equation}
\begin{equation} \label{sys2du}
\begin{aligned}
&\dot{x}^{u,j} = f^{u,j}(x^{u,j},-\omega_j), \\
&-d^u_j = g^{u,j}(x^{u,j},-\omega_j),
\end{aligned} \hspace{2em}j \in N
\end{equation}
\end{subequations}
\n{where the input $\zeta_j$ is defined as $\zeta_j = [ -\omega_j \; p_j^c]^T$ with $p_j^c$ representing the deviations of a power command signal from its nominal value. {\color{black}Notice that {\color{black}in the case of} uncontrollable demand,} the input is given in terms of the local frequency deviation $\omega_j$ only, and is decoupled from the power command signal as expected.}


For notational convenience, we collect the variables in~\eqref{sys2} into the vectors $x^M = [x^{M,j}]_{j \in G}$, $x^c = [x^{c,j}]_{j \in N}$, and $x^u = [x^{u,j}]_{j \in N}$.
These quantities represent the internal states of the dynamical systems used to update the {\color{black}outputs} $p^M_j$, $d^c_j$, and $d^u_j$.


In terms of the outputs from~\eqref{sys2}, it will be useful to consider the net supply variables $s$, {\color{black} defined as
{\color{black}
\begin{equation} \label{ssys} 
s_j = p^M_j - d^c_j , \; j \in G,
\qquad s_j = -d^c_j , \; j \in L.
\end{equation}
}
The variables defined in~\eqref{ssys} evolve according to the dynamics described in \eqref{sys2p} - \eqref{sys2dc}. Therefore, $s_j$ {\color{black} are outputs} from these combined {\color{black}controlled dynamical} systems with {\color{black}inputs $\zeta_j$.}}

\vspace{-2mm}
\subsection{Power Command Dynamics}\label{Power_Command_Dynamics}

{\color{black}
We {\color{black}consider a communication network described by} a connected graph {\color{black}($N,\tilde E$)}, where {\color{black}$\tilde E$} represents the set of communication lines among the buses, {\color{black}i.e., $(i,j) \in \tilde E$ if buses} $i$ and $j$ communicate. {\color{black}Note that $\tilde E$} can be different from the set of flow lines $E$.}  {\color{black}We will study the behavior of the system \eqref{sys1}--\eqref{sys2} under the following {\color{black}dynamics} for the power command signal $p^c_j$
which has been used in literature (e.g. \cite{P30_Low,P6_Low}),}
{\color{black}
\begin{subequations}\label{sys_power_command}
\begin{equation}\label{sys_power_command_a}
\gamma_{ij} \dot{\psi}_{ij} = p_i^c - p_j^c , \; (i,j) \in {\color{black}\tilde E}
\end{equation}
\begin{equation}\label{sys_power_command_b}
{\color{black}\gamma_j \dot{p}_j^c = -(s_j - p_j^L)  - \sum_{k:j\rightarrow k} \psi_{jk} + \sum_{i:i\rightarrow j} \psi_{ij}, \; j \in N}
\end{equation}
\end{subequations}
where} {\color{black}$\gamma_j$ and $\gamma_{ij}$ are positive constants, and the variable $\psi_{ij}$ represents the difference in the integrals between the power commands of communicating buses $i$ and $j$}. It should be noted that $p_i^c$ and $p_j^c$ are {\color{black}variables} shared between {\color{black}communicating} buses $i$ and~$j$.

 Although the dynamics in \eqref{sys_power_command} do not directly integrate frequency, {\color{black} we will see later that under a weak condition on the {\color{black}steady state {\color{black}behavior} of $d^u,$} they guarantee convergence to {\color{black} the nominal frequency}} for a broad class of supply dynamics.
{\color{black}The dynamics in \eqref{sys_power_command}}, often referred as `virtual swing equations', are frequently {\color{black}used in the literature\footnote{\color{black}In this paper we use for simplicity a single communicating variable. It should be noted that more advanced communication structures (e.g. \cite{P30_Low}) can allow additional constraints to be satisfied in the optimization problem posed.}} as they achieve both the {\color{black}synchronization} of {\color{black}the communicated variable  {\color{black}$p^c$}}, something {\color{black}that can be exploited to guarantee optimality of the equilibrium point {\color{black}reached,
and} {\color{black}also} the convergence of frequency to its nominal value.

\subsection{Optimal supply and load control} \label{sec:optim}

We aim to study how generation and controllable demand should be adjusted in order to meet the step change in frequency independent \n{demand} and simultaneously {\color{black}minimize} the cost that comes from the deviation {\color{black}in the power} generated and the disutility of loads. We now introduce an optimization problem, which we call the optimal supply and load control problem (OSLC), {\color{black}that} can be used to achieve this goal.

A cost $C_j(p_j^M)$ is supposed to be incurred when~generation output at bus $j$ is changed by $p_j^M$ from its nominal value. Similarly, a cost of {\color{black}$C_{dj} (d^c_{j})$ is incurred for a change of $d^c_{j}$} in controllable demand.
The total cost within OSLC is the sum of the above costs. The problem is to find the vectors $p^M$ and $d^c$ that minimize this total cost and simultaneously achieve power balance, while satisfying physical saturation constraints. {\color{black}More precisely, the  following optimization problem is considered}

\begin{equation}
\begin{aligned}
&\hspace{-2em}\underline{\text{OSLC:}} \\
&\min_{p^M,d^c} \sum\limits_{j\in G} C_{j} (p_{j}^M) + \sum\limits_{j\in N}  C_{dj} (d^c_{j}),  \hspace{-1.5em}\\
&\text{subject to } \sum\limits_{j\in  G} p_j^M = \sum\limits_{j\in  N} (d^c_j + p_j^L), \\
& p^{M,min}_j \leq p^M_j \leq p^{M,max}_j \,,\, \forall j \in G, \\
& d^{c,min}_j \leq d^c_j \leq d^{c,max}_j \,,\, \forall j \in N,
 \label{Problem_To_Min}
\end{aligned}
\end{equation}
where $p^{M,min}_j,p^{M,max}_j, d^{c,min}_j$, and $d^{c,max}_j$ are bounds for the minimum and maximum values for generation and controllable demand deviations, respectively, at bus~$j$.
The equality constraint in~\eqref{Problem_To_Min} requires all the \n{additional} frequency-independent loads to be matched by the total deviation in generation and controllable demand.
This {\color{black} ensures that when system \eqref{sys1} is at equilibrium and {\color{black}a mild} condition described in Assumption \ref{assum2.4} below holds, the frequency will be at its nominal value.}

{\color{black} Within the paper we aim to specify properties on the control dynamics of $p^M$ and $d^c$, described in \eqref{sys2p}--\eqref{sys2dc},
that ensure that those quantities 
converge to values at which optimality can be guaranteed for \eqref{Problem_To_Min}.}




{\color{black}The assumption below
allows the use of the KKT conditions to prove the optimality result in Theorem~\ref{optthm} in Section \ref{Main_Results}.

\begin{assumption} \label{assum5}
The cost functions $C_{j}$ and $C_{dj}$ are continuously differentiable and strictly convex.
\end{assumption}
}

\subsection{Equilibrium analysis}

We now {\color{black}describe} what is meant by an {\color{black}equilibrium}
of the interconnected system~\eqref{sys1}--\eqref{sys_power_command}.

\begin{definition} \label{eqbrdef}
{\color{black}The point $\beta^* = (\eta^*,\psi^*, \omega^*,x^{M,*}, x^{c,*},$ $ x^{u,*}, p^{c,*})$ defines} an equilibrium of the system~\eqref{sys1}--\eqref{sys_power_command} if all time derivatives of ~\eqref{sys1}--\eqref{sys_power_command} are equal to zero {\color{black}at this point}.
%

\end{definition}

It should be noted that the static input-output maps $k_{p^M_j}$, $k_{d^c_j}$, and $k_{d^u_j}$, as defined in \eqref{iochar}, completely {\color{black}characterize} the equilibrium {\color{black} behavior} of \eqref{sys2}.
{\color{black}In our analysis, we shall consider conditions on these characteristic maps relating input {\color{black}$\zeta_j = [ -\omega_j \; p_j^c]^T$} and generation/demand such that their equilibrium values {\color{black}are optimal for \eqref{Problem_To_Min}, thus making sure that frequency} will be at its nominal value at steady state}.

Throughout the paper, it is \n{assumed} that there exists some equilibrium of~\eqref{sys1}--\eqref{sys_power_command} as defined in Definition~\ref{eqbrdef}. Any such equilibrium is denoted by $\beta^* = (\eta^*,\psi^*, \omega^*, x^{M,*}, x^{c,*}, x^{u,*}, p^{c,*})$. {\color{black}Furthermore, we use $(p^*, p^{M,*}, d^{c,*}, d^{u,*}, \zeta^*, s^{*})$ to represent the equilibrium values of respective quantities
{\color{black}
in~\eqref{sys1}--\eqref{sys_power_command}.
}

{\color{black}The power angle differences at the considered equilibrium are assumed to satisfy the following security constraint.}

\begin{assumption} \label{assum1}
$| \eta^*_{ij} | < \tfrac{\pi}{2}$ for all $(i,j) \in E$.
\end{assumption}

Moreover, the following {\color{black}assumption is related} with the steady state values of variable $d^u$, describing uncontrollable demand and generation damping. {\color{black} It is a mild condition associated with having negative feedback from $d^u$ to frequency.}

\begin{assumption}\label{assum2.4}
 For \n{each} $j \in N$, the functions $k_{d^u_j}$ relating the steady state values of frequency and uncontrollable loads satisfy {\color{black}$\bar u_j k_{d^u_j}(\bar u_j) > 0$ for all $\bar u_j \in \mathbb{R}-\{0\}$}.
\end{assumption}}

{\color{black} Although not required for stability, Assumption \ref{assum2.4} guarantees that {\color{black}the} frequency will be equal to its nominal value at equilibrium, i.e. $\omega^{*} = \vect{0}_{|N|}$, as stated in \n{the following lemma}, {\color{black}proved} in {\color{black}Appendix A.}

 \begin{lemma}\label{eqbr_lemma}
 Let Assumption \ref{assum2.4} hold. Then, any equilibrium point $\beta^*$ given by Definition \ref{eqbrdef} satisfies $\omega^{*} = \vect{0}_{|N|}$.
 \end{lemma}

The stability and optimality properties of such equilibria will be studied in the following sections.}

\subsection{Additional conditions}
{\color{black}
Due to the fact that the frequency at the load buses is related with the system {\color{black}states} by means of algebraic equations, additional conditions are needed for the system \eqref{sys1}--\eqref{sys2} to be well-defined.}
\n{\color{black}We use below the}
vector notation $\omega^G = [\omega_j]_{j \in G}$ and $\omega^L = [\omega_j]_{j \in L}$.}
\begin{assumption} \label{assum4}
{\color{black}There exists an open neighborhood $T$ of $(\eta^*, \omega^{G,*}, x^{c,*},x^{u,*}, p^{c,*})$
and a locally Lipschitz map $f^L$ such that when  $(\eta, \omega^G, x^c, x^u, p^c) \in T$,
$\omega^L = f^L(\eta, \omega^G, x^c, x^u, p^c)$.}
%
\end{assumption}

\begin{remark}
Assumption \ref{assum4} is a technical assumption that is
required in order for the system \eqref{sys1}--\eqref{sys2} to have a locally
well-defined state space realization.
{\color{black}It can often be easily verified by means of the implicit function theorem \cite{Hill_Mareel}}.
Without Assumption \ref{assum4}, stability could be \n{studied} through more technical  {\color{black}approaches} such as the
singular perturbation analysis discussed in \cite[Section 6.4]{sastry}.
\end{remark}
}

{\color{black}\section{Dissipativity conditions on generation and demand dynamics} \label{Assumptions}}

{\color{black} Before we state our main results in Section \ref{Main_Results}, it would be useful to}
provide a {\color{black}dissipativity definition, based on \cite{Dissipativity_book},} for systems {\color{black}of the form} \eqref{dynsys}. This notion will be used to {\color{black} formulate appropriate decentralized conditions on the uncontrollable demand and power supply {\color{black}dynamics} \eqref{sys2du}, \eqref{ssys}}.

\begin{definition}\label{Dissipativity_Definition}
The system \eqref{dynsys} is said to be locally dissipative about the constant input values $\bar{u}$ and {\color{black}corresponding equilibrium} state values $\bar{x}$,  with supply rate function $W: \mathbb{R}^{n+k} \to \mathbb{R}$, if there exist open {\color{black}neighborhoods} $U$
of $\bar{u}$ and $X$ of $\bar{x}$, and  a continuously differentiable, {\color{black}positive definite function {\color{black}$V: \mathbb{R}^{m} \to \mathbb{R}$ (called the storage function), with a strict local minimum at $x=\bar{x}$, such that for all $u \in U$ and all $x \in X$,}}
\begin{equation}\label{dissipativity_condition}
\dot{V}(x) \leq W(u,y).
\end{equation}
\end{definition}

We now {\color{black}assume} that {\color{black}the systems with} input~$\zeta_j = [ -\omega_j \;\;  p^c_j]^T$  {\color{black}and output the} power supply variables and uncontrollable loads satisfy
{\color{black}the following local dissipativity condition.}

\begin{assumption}\label{assum2}
{\color{black}The systems with inputs $\zeta_j = [ -\omega_j \; p^c_j]^T$ and outputs $y_j = [s_j$ $-d^u_j]^T$ described in \eqref{ssys} and \eqref{sys2du}
satisfy a dissipativity condition
about constant input values $\zeta_j^*$ and corresponding equilibrium state values $(x^{M,j,*},x^{c,j,*}, x^{u,j,*})$ in the sense \n{of} Definition \ref{Dissipativity_Definition}, with supply rate functions
\begin{multline}\label{supply_rate}
W_j(\zeta_j,y_j) = [(s_j - s^{*}_j) \;\;\;  (-d^u_j - (-d^{u,*}_j))] \begin{bmatrix}
1 & 1\\
1 & 0
\end{bmatrix} (\zeta_j - \zeta^*_j) \\
-\phi_j(\zeta_j - \zeta^*_j),\, j \in N.
\end{multline}
Furthermore, one of the following two properties holds,
\begin{enumerate}[(a)]
\item The function $\phi_j$ is positive definite.
\item The function $\phi_j$ is positive semidefinite and positive definite with respect to
$\omega_j$. Also
     when $\omega_j$, $s_j$ are constant for all times then $p_j^c$ cannot be a nontrivial sinusoid\footnote{{\color{black}By nontrivial sinusoid, we mean functions of the form $\sum_j A_j \sin(\omega_j t+\phi_j)$ {\color{black}that are not equal to a constant.}}}.
\end{enumerate}
}
\end{assumption}
{\color{black}
We shall refer to {\color{black}Assumption} \ref{assum2} when condition (a) holds for $\phi_j$ as {\color{black}Assumption} \ref{assum2}(a) (respectively Assumption \ref{assum2}(b) when (b) holds).}

\begin{remark}
{\color{black}Assumption~\ref{assum2} is a decentralized condition that allows to incorporate a broad class of generation and load dynamics, including various examples that have been used in the literature (these will be discussed in Section~\ref{Discussion}).
{\color{black}Furthermore, for linear systems Assumption~\ref{assum2}} can be formulated as the feasibility problem of a corresponding LMI (linear matrix inequality) \cite{Dissipativity_book_2}, and it can therefore be verified by means of computationally efficient methods. }
\end{remark}

\begin{remark}\label{rem:imzeros}
Condition (b) in Assumption \ref{assum2} is a relaxation of condition (a) whereby $\phi$ is not required to be positive definite. This permits the inclusion of a broader class of dynamics from $p^c_j$ to $s_j$ as {\color{black}it} will be discussed in Section~\ref{Discussion}.  However, it requires that {\color{black}the} {\color{black}power command $p^c$} cannot be a sinusoid if both $s_j$ and $\omega_j$ are {\color{black}constant}. This {\color{black}additional {\color{black}condition} is necessary} as the dynamics in \eqref{sys_power_command} allow $p^c_j$ {\color{black}to be a sinusoid} when $s_j$ is {\color{black}constant}. For linear systems, {\color{black}this condition is implied by the} rather mild assumption that no imaginary axis zeros are present in the transfer function from $p^c_j$ to $s_j$.
\end{remark}

\begin{remark}
{\color{black}Further intuition on the dissipativity condition in Assumption \ref{assum2}  will be provided in Section \ref{System_Representation} and Appendix~B. In particular, it will be shown when $\phi_j=0$ that this is a decentralized condition that is necessary and sufficient} for the {\color{black}passivity}
of {\color{black}an appropriately defined multivariable system quantifying aggregate dynamics at each bus.}
\end{remark}

\section{Main Results}\label{Main_Results}

In this section  we
state our main results, with their proofs provided in {\color{black}Appendix A.}


{\color{black}


Our first result provides {\color{black}conditions} for the equilibrium points to be solutions\footnote{\color{black}
Note that an equilibrium point is a solution to the OSLC problem when at that point the variables that appear in \eqref{Problem_To_Min} are solutions to the problem.} to the OSLC problem \eqref{Problem_To_Min}.

\begin{theorem} \label{optthm}
Suppose that Assumption \ref{assum5} is satisfied. If the control dynamics in \eqref{sys2p} and \eqref{sys2dc} are chosen such that
\begin{equation}\label{contspec}
\begin{aligned}
k_{p^M_j}(\zeta_j) &= [(C_j')^{-1}(f(\zeta_j)]^{p^{M,max}_j}_{p^{M,min}_j}\\
 k_{d^c_j}(\zeta_j) &= [(C_{dj}')^{-1}(-f(\zeta_j)]^{d^{c,max}_j}_{d^{c,min}_j}
\end{aligned}
\end{equation}
holds for some  surjective function $f:\mathbb{R}^2\rightarrow \mathbb{R}$, {\color{black}which is {\color{black}strictly} increasing with respect to $p^c_j$,}
then the equilibrium values $p^{M,*}$ and $d^{c,*}$
are optimal for the OSLC problem \eqref{Problem_To_Min}.
\end{theorem}

Our second result shows that the set of equilibria for the system described by \eqref{sys1}--\eqref{sys_power_command} for which Assumptions \ref{assum5} - \ref{assum2} are satisfied is asymptotically attracting{\color{black}, the equilibria are global minima of the OSLC} problem \eqref{Problem_To_Min} and, as shown in Lemma~\ref{eqbr_lemma}, {\color{black}satisfy} $\omega^{*} = \vect{0}_{|N|}$.

\begin{theorem} \label{convthm}
Consider equilibria of~\eqref{sys1}--\eqref{sys_power_command} with respect to which Assumptions~\ref{assum5}--\ref{assum2} are all satisfied. If the control dynamics in~\eqref{sys2p} and~\eqref{sys2dc} are chosen such that \eqref{contspec} holds,
%
then there exists an open neighborhood of initial conditions about any such equilibrium such that the solutions of~\eqref{sys1}--\eqref{sys_power_command} are guaranteed to converge to a
{\color{black}set of equilibria that solve} the OSLC problem~\eqref{Problem_To_Min} with $\omega^* = \vect{0}_{|N|}$.
\end{theorem}
}

\section{Discussion}\label{Discussion}

{\color{black} In this section we discuss examples that fit within the framework presented in the paper, and also describe how the dissipativity condition of Assumption \ref{assum2} can be verified for linear systems via a {\color{black}linear matrix inequality.}
}

We start by giving various examples of power supply dynamics that have been {\color{black}used in the literature} that satisfy \n{our proposed} dissipativity condition in Assumption \ref{assum2}.
Consider the {\color{black}load models} used in  {\color{Black}\cite{P30_Low}, \cite{P21}, and \cite{P6_Low},} where {\color{black}the power} supply is {\color{black}a static function of {\color{black}$\omega_j$} and $p^c_j$},  {\color{black}
\begin{equation}\label{example_static}
s_j = (C_j')^{-1}(p^c_j - \omega_j), \,\, j \in N,
\end{equation}
where $C_j$ is some convex cost function, and {\color{black}generation damping/uncontrollable demand} is {\color{black}given} by $d^u_j = \lambda_j \omega_j, \lambda_j > 0$.
It is easy to show that Assumption \ref{assum2}(a) holds for these widely used schemes}.

Furthermore, Assumption \ref{assum2}(b) is satisfied
{\color{black}when first order generation dynamics are used such as
 \begin{equation}\label{sys_example_Low}
{\color{black}\dot{s}_j = -\mu_j(C_j'(s_j) -{\color{black}(p^c_j - \omega_j)})}
\end{equation}
with $d^u_j = \lambda_j \omega_j$ {\color{black}and} $\lambda_j,\mu_j > 0$. Such first order models have often been used in the literature as in \cite{P1_Low}.}


{\color{black}A significant aspect of the framework presented in this paper is that it also allows higher} order dynamics for {\color{black}the} power supply {\color{black} to be {\color{black}incorporated.}
{\color{black} {\color{black}As an example, we consider the following second-order model,}
{\color{black}
\vspace{-1mm}
\begin{equation} \label{sys_example}
\begin{aligned}
& \dot{\alpha}_j = -\frac{1}{\tau_{a,j}}(\alpha_j - K_j(p_j^c - \omega_j)), \\
& \dot{z}_j = -\frac{1}{\tau_{b,j}}(z_j - \alpha_j), \\
& {\color{black}s_j - d^u_j = z_j -\lambda_j\omega_j+\lambda_j^{PC} p^c_j,}
\end{aligned}
\end{equation}
where $\alpha_j$, $z_j$ are states and $\tau_{a,j},\tau_{b,j}>0$ time constants associated with the turbine-governor dynamics,  $\lambda_j>0$ is a damping coefficient\footnote{\color{black} Note that the term $\lambda_j\omega_j$ {\color{black}can be incorporated
in $s_j$ or~$d^u_j$.}}, {\color{black}constant $K_j>0$ determines the strength of the feedback gain, and
the term $\lambda_j^{PC} p^c_j$ {\color{black}represents static dependence on power command due to either generation or controllable loads\footnote{\color{black} It should be noted that the term $d^u_j$ can also include controllable demand and generation that depend on frequency only (i.e. not on power command). Therefore, $d^u_j$ can be perceived to contain all frequency dependent {\color{black}terms} that return to their nominal value at steady state and therefore do not contribute to secondary frequency control.}.}
 It can be shown that Assumption~\ref{assum2} is satisfied for all $\tau_{a,j},\tau_{b,j}>0$
 when\footnote{\color{black}A second order model was studied for a related problem in \cite{P57_Low}, with the stability condition requiring, roughly speaking, that the gain of the system is less than the damping provided by the loads. The LMI approach described in this section allows such conditions to {\color{black}be relaxed.}}
 $K_j<8\lambda_j^{PC}$ and $\lambda^{PC}_j \leq \lambda_j$}.}

{\color{black}
\n{Another} feature of Assumption \ref{assum2} is that it can be efficiently verified for a general linear system by means of an LMI, i.e. a computationally efficient convex problem. In particular, it can be shown \cite{Dissipativity_book_2} that if the system in Assumption \ref{assum2} is linear with a minimal state space realization
\begin{equation}\label{sspace_eq}
\begin{aligned}
&\dot{x} = Ax + B\tilde u, \\
&\tilde y = Cx + D\tilde u,
\end{aligned}
\end{equation}
where $\tilde u= \zeta-\zeta^\ast$ and $\tilde y= y-y^\ast$, {\color{black}and $\phi_j$} is chosen as a quadratic function $\phi_j=\epsilon_1(\omega_j-\omega_j^\ast)^2+\epsilon_2(p^c_j-p_j^{c,\ast})^2$ with\footnote{\color{black}We could also have $\epsilon_2=0$ if \eqref{sspace_eq} has no zeros on the {\color{black}imaginary} axis, as stated in condition (b) for $\phi$ in Assumption~\ref{assum2}, and Remark \ref{rem:imzeros}.}  {\color{black}$\epsilon_1,\epsilon_2>0$}
 then the dissipativity condition in {\color{black}Assumption~\ref{assum2}} is satisfied if and only if there exists $P=P^T\geq0$ such that
\begin{equation}\label{LMI}
\begin{bmatrix}
A^T P + PA  & PB \\
B^TP & 0
\end{bmatrix}
-
\begin{bmatrix}
C  & D \\
0 & I
\end{bmatrix}^T
Q
\begin{bmatrix}
C  & D \\
0 & I
\end{bmatrix}
\leq 0,
\end{equation}
where the matrix $Q$ is given by
 \begin{equation*} 
 Q=
  		\begin{bmatrix}
          0  & M \\
          M & K
          \end{bmatrix},\
  M={\color{black} \frac{1}{2}}
          \begin{bmatrix}
          1 & 1 \\
          1 & 0
          \end{bmatrix},\
          K=
		\begin{bmatrix}
          -\epsilon_1 & 0 \\
          0 & -\epsilon_2
          \end{bmatrix}{\color{Black}.}
 \end{equation*}

%

{\color{black} This approach could also be exploited to form various convex {\color{black}optimization} problems that could facilitate design.  For example, one could obtain the minimum {\color{black}damping such that}
Assumption \ref{assum2} is satisfied {\color{black}at a bus}}.}


{\color{black}
To further demonstrate the applicability of our approach we consider a fifth order model for turbine governor dynamics provided by the Power System Toolbox \cite{19-25}.
The dynamics are described by the following transfer function relating {\color{black}{\color{black}the mechanical power} supply\footnote{\color{black}Note that {\color{black}$\hat{s}_j$} denotes the Laplace transform of {\color{black}$s_j$}.} $\hat{s}_j$  with the negative frequency deviation $-\hat{\omega}_j$,}
\[
G_j(s)=K_j\frac{1}{(1+sT_{s,j})}\frac{(1+sT_{3,j})}{(1+sT_{c,j})}\frac{(1+sT_{4,j})}{(1+sT_{5,j})},
\]
where $K_j$ and $T_{s,j}, T_{3,j}, T_{c,j}, T_{4,j}, T_{5,j}$ are the droop coefficient and time-constants respectively. Realistic values for these models are provided by the toolbox for the NPCC network, with turbine governor dynamics implemented on 22 buses.
The corresponding buses also have appropriate frequency damping $\lambda_j$.
{\color{black}We examined the effect of incorporating a power command input signal in the above dynamics by considering the {\color{black}{\color{black}supply dynamics}
\[
\hat{s}_j -\hat{d}^u_j = (G_j + \lambda_j) (-\hat{\omega}_j) + (G_j + \lambda^{PC}_j) (\hat{p}^c_j), j \in N
\]
where $\lambda^{PC}_j>0, j \in N$} is a coefficient representing the static dependence on power command.
 {\color{black}For} appropriate values of $\lambda^{PC}_j$, the condition {\color{black}in Assumption \ref{assum2}} was satisfied for 20 out of the 22 buses, while for the remaining 2 buses the damping coefficients $\lambda_j$ needed to be increased by $37\%$ and $28\%$ respectively.}
 {\color{black}
 Furthermore, {\color{black} filtering  the power command signal with appropriate {\color{black}compensators,}}
allowed a significant decrease in the required }\color{black}value}
for $\lambda^{PC}_j$.
 {\color{black} Power command {\color{black}and frequency} compensation may also be used with alternative objectives, such as to {\color{black}improve} the stability margins {\color{black}and} system performance.}

{\color{black}The fact that our condition {\color{black}is satisfied at all but two buses\footnote{\color{black}Note that this is satisfied at all buses with appropriate increase in damping.}, {\color{black}demonstrates that {\color{black}it} is not conservative in existing implementations. Note also that a main feature of this condition is the fact that {\color{black}it is decentralized,} involving only local bus dynamics, which can be important in practical implementations.}}

}

}

{\color{black}
\subsection{System Representation}\label{System_Representation}

It is useful and intuitive to note that {\color{black}the system \eqref{sys1} - \eqref{sys_power_command} considered in the paper} can be represented by a negative feedback interconnection of {\color{black}systems $I$ and {\color{black}$B=\bigoplus_{j=1}^{|N|} B_j$},} containing all interconnection and bus dynamics respectively. {\color{black}More precisely,} $I$ and $B$ have respective inputs $u^I$ and $u^B$, and {\color{black}respective} outputs $u^B$ and $-u^I$, defined as
\begin{equation*}
u^I \hspace{-0.5mm} = \hspace{-0.5mm}
	\begin{bmatrix}
    \omega_1  \\
    -p^c_1 \\
    \dots \\
    \omega_{|N|}  \\
    -p^c_{|N|}
	\end{bmatrix}, \
 u^B\hspace{-0.5mm} = \hspace{-0.5mm}
    \begin{bmatrix}
    \sum_{k:1\rightarrow k} p_{1k} - \sum_{i:i\rightarrow 1} p_{i1}   \\
    \sum_{i:i\rightarrow 1} \psi_{i1} - \sum_{k:1\rightarrow k} \psi_{1k} \\
    \dots \\
    \sum_{k:|N|\rightarrow k} p_{|N|k} - \sum_{i:i\rightarrow |N|} p_{i|N|} \\
    \sum_{i:i\rightarrow |N|} \psi_{i|N|} - \sum_{k:|N|\rightarrow k} \psi_{|N|k}
	\end{bmatrix}.
\end{equation*}}
%
%
%
{\color{black}
The subsystems $B_j$ represent the dynamics at bus $j$ and have  inputs $[u^B_{2j-1} \, u^B_{2j}]^T$ and outputs $[-u^I_{2j-1}$ $-u^I_{2j}]^T$.

{\color{black}
{\color{black}It can easily be shown that System $I$ is locally passive\footnote{\color{black}By a locally passive system we refer to a system satisfying the dissipativity condition in Definition \ref{Dissipativity_Definition} with the supply rate being $W(u,y)=(u-u^*)^T(y-y^*)$.}}.
{\color{black} The following theorem shows that Assumption \ref{assum2}  with  $\phi=0$
is sufficient for the passivity of each individual subsystem $B_j$.
\begin{theorem}\label{suff_thm}
Consider the system described by \eqref{sys1} - \eqref{sys_power_command} and its representation by systems $I$ and $B$, defined in section \ref{System_Representation}.
Then, the dissipativity condition in Assumption \ref{assum2}  with $\phi=0$ is sufficient for the passivity of subsystems $B_j, j\in N$ about the equilibrium point considered in Assumption \ref{assum2}.
\end{theorem}}

\begin{remark}
The significance of the interpretation discussed in this section is that the passivity property of system $I$, in conjunction with the fact that $B=\bigoplus_{j=1}^{|N|} B_j$,
implies that stability of the network is guaranteed if the subsystems $B_j$ are passive (with appropriate strictness as quantified within the paper). In particular, stability is guaranteed in a decentralized way without requiring information about the rest of the network at each individual bus, which \n{is advantageous} in highly distributed schemes where a "plug and play" capability is needed. It should be noted that the subsystems $B_j$ are multivariable systems quantifying the aggregate bus dynamics associated with both power generation and the communicated signal~$p^c$.
\end{remark}

{\color{black}
\begin{remark}
It is shown in Appendix B that Assumption \ref{assum2} is {\color{black}also necessary}
for systems~$B_j$ to be passive, \n{for general} affine nonlinear dynamics.
Hence,  Assumption \ref{assum2} introduces no additional conservatism {\color{black}in this property} for a large class of {\color{black}nonlinear systems}.
\end{remark}
}


{\color{black}
\subsection{Observing uncontrollable frequency independent demand}\label{sec: Observer}

{\color{black}The power command dynamics  in \eqref{sys_power_command} involve the uncontrollable frequency independent demand $p^L$.
We discuss in this section that the inclusion of appropriate observer dynamics for $p^L$ allows convergence to optimality to be achieved when $p^L$ is not directly known.}

A way to obtain $p^L$, could be by re-arranging equations \eqref{sys1b}--\eqref{sys1c}.
This approach  would require knowledge of power supply and power transfers in load buses, which is realistic. However, knowledge of the frequency derivative would also be required for its estimation at generation buses, which might be difficult to obtain \n{in noisy environments}.

{\color{black}We therefore consider instead observer dynamics\footnote{\color{black}See also the use of observer dynamics in \cite{nima_depersis} as a means of counteracting agent dishonesty.} for $p^L_j$ that are incorporated within the power command dynamics. In particular the following dynamics are considered}

\begin{subequations}\label{sys_observer}
\begin{equation}\label{sys_observer_a}
\gamma_{ij} \dot{\psi}_{ij} = p_i^c - p_j^c , \; (i,j) \in \tilde E,
\end{equation}
\begin{equation}\label{sys_observer_b}
\gamma_j \dot{p}_j^c = -(s_j - \chi_j)  - \sum_{k:j\rightarrow k} \psi_{jk} + \sum_{i:i\rightarrow j} \psi_{ij}, \; j \in N,
\end{equation}
\begin{equation}\label{sys_observer_c}
\tau_{\chi,j}\dot{\chi}_j = b_j - \omega_j - p^c_j - \chi_j, \; j \in G,
\end{equation}
\begin{equation}\label{sys_observer_d}
M_j \dot{b}_j = - \chi_j + s_j - d^u_j - \sum_{k:j\rightarrow k} p_{jk} + \sum_{i:i\rightarrow j} p_{ij}, \; j\in G,
\end{equation}
\begin{equation}\label{sys_observer_e}
 0 = - \chi_j + s_j - d^u_j - \sum_{k:j\rightarrow k} p_{jk} + \sum_{i:i\rightarrow j} p_{ij}, \; j\in L,
\end{equation}
\end{subequations}
where $\tau_{\chi,j}$ are positive time constants and $b_j$ and $\chi_j$ {\color{black}are}
auxiliary variables associated with the observer.



{\color{black}
The equilibria of the system \eqref{sys1} -- \eqref{ssys}, \eqref{sys_observer} are defined in a similar way to Definition \ref{eqbrdef} and it is assumed that at least one such equilibrium exists.
 Note that the existence of an equilibrium of \eqref{sys1} - \eqref{sys_power_command} implies the existence of an equilibrium of ~\eqref{sys1}--\eqref{ssys}, \eqref{sys_observer}.}

\n{We now {\color{black}provide a result analogous to Lemma \ref{eqbr_lemma} in} the case where the observer dynamics {\color{black}are} included.} {\color{black} Lemma \ref{eqbr_lemma_observer} is proven in  Appendix A.}

  \begin{lemma}\label{eqbr_lemma_observer}
 Let Assumption \ref{assum2.4} hold. Then, any equilibrium point $(\eta^*$,$\psi^*,\omega^*, x^{M,*}, x^{c,*}$, $x^{u,*},p^{c,*}$,$b^*$,$\chi^*)$ of the system~\eqref{sys1} -- \eqref{ssys}, \eqref{sys_observer} satisfies $\omega^{*} = \vect{0}_{|N|}$.
 \end{lemma}

{\color{black}
\begin{remark}
The dynamics in \eqref{sys_observer} eliminate the requirement to explicitly know $p^L$ within the power command dynamics by adding an observer that mimics the swing equation, described by \eqref{sys_observer_c}--\eqref{sys_observer_e}.
The dynamics in \eqref{sys_observer_d}--\eqref{sys_observer_e} ensure that the variable $\chi_j$ 
{\color{black}is equal at steady state to}
the value $\chi^*_j = s^*_j - d^{u,*}_j - \sum_{k:j\rightarrow k} p^*_{jk} + \sum_{i:i\rightarrow j} p^*_{ij} = p^L_j$ for $j \in N$, with the second part of the equality coming from \eqref{sys1b}--\eqref{sys1c} at equilibrium.
 As shown in {\color{black}Lemma \ref{eqbr_lemma_observer},} such equilibrium guarantees that the steady state value of {\color{black}the} frequency will be equal to the {\color{black}nominal {one}.}
\end{remark}
}

{\color{black}
The following proposition, proved in Appendix A, shows that the
 set of equilibria for the system described by \eqref{sys1} -- \eqref{ssys}, \eqref{sys_observer} for which Assumptions \ref{assum5} - \ref{assum2} are satisfied is asymptotically attracting and {\color{black}that} these equilibria are also solutions to the OSLC problem \eqref{Problem_To_Min}.
\begin{proposition}\label{convthm_observer}
Consider equilibria of~\eqref{sys1} -- \eqref{ssys}, \eqref{sys_observer}  with respect to which Assumptions~\ref{assum5}--\ref{assum2} are all satisfied. If the control dynamics in~\eqref{sys2p} and~\eqref{sys2dc} are chosen such that \eqref{contspec} is satisfied
then there exists an open neighborhood of initial conditions about any such equilibrium such that the solutions of~\eqref{sys1} -- \eqref{ssys}, \eqref{sys_observer}  are guaranteed to converge to a global minimum of the OSLC problem~\eqref{Problem_To_Min} with $\omega^* = \vect{0}_{|N|}$.
\end{proposition}
}

{\color{black}
\begin{remark}
{\color{black}Note that in some cases there could be uncertainty in the knowledge of the $d^u$ dynamics. This does not affect the optimality of the equilibrium points since at equilibrium we {\color{black}have $d^u=\vect{0}_{|N|}$}. Numerical simulations with realistic data have demonstrated that network stability is also robust to variations in the $d^u$ model used in  \eqref{sys_observer_d}--\eqref{sys_observer_e}.  }

\end{remark}
}

\vspace{-1mm}
\section{{\color{black}Simulation on the NPCC 140-bus} system} \label{Simulation}

In this section we use the {\color{black}Northeast Power Coordinating
Council (NPCC) 140-bus interconnection system,} simulated using the Power System Toolbox~\cite{19-25}, in order to illustrate our results. This model is more detailed and realistic than our analytical one, including line resistances, a DC12 exciter model{\color{black}, a subtransient reactance generator model, and {\color{black}higher} order turbine governor models}\footnote{The details of the simulation models can be found in the Power System Toolbox data file {\color{black}datanp48}.}.

\begin{figure}[t]
\centering
\includegraphics[trim = 1mm 0mm 10mm 0mm, height = 2.2in,width=3.115in,clip=true]{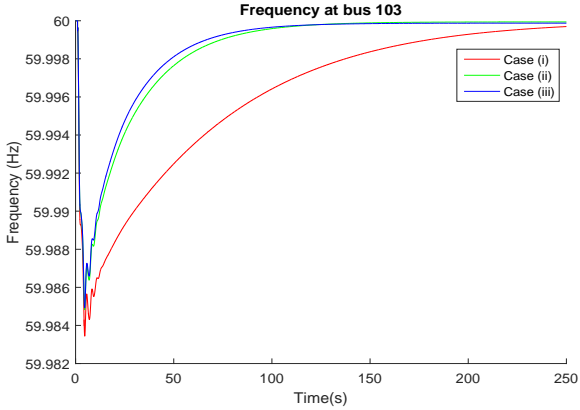}
\vspace{-1mm}
\caption{\color{black}Frequency at bus 103 with: i) 10 generators, ii) 10 generators and 20 controllable loads, iii) {\color{black}15} generators and 20 controllable loads, contributing to secondary frequency control.} \vspace{-4.5mm}
\label{Frequency}
\end{figure}


\begin{figure}[t]
\centering
\includegraphics[trim = 1mm 0mm 10mm 0mm, height = 2.2in,width=3.115in,clip=true]{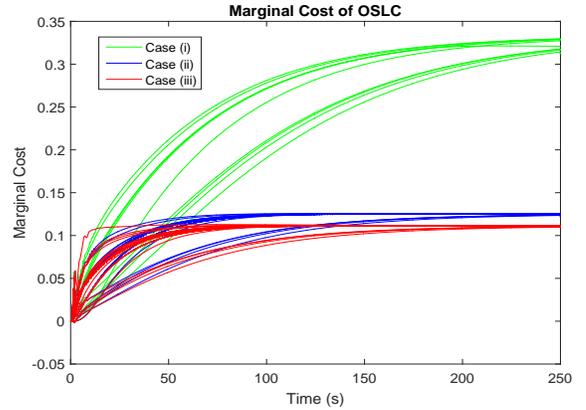}
\vspace{-1mm}
\caption{\color{black}Marginal costs for controllable loads and generators with non-equal cost coefficients for the three {\color{black}test} cases. } \vspace{-5mm}
\label{Marginal}
\end{figure}

The test system consists of {\color{black}93 load} buses serving different types of loads including constant active and reactive loads and {\color{black}47 generation} buses. The overall system has a total real power of {\color{black}28.55GW.} {\color{black}For our simulation, we added three loads on units 2, 9, and 17,} each having a step increase of magnitude $1$ p.u. (base 100MVA) at $t=1$ second.

{\color{black}
Controllable demand was considered within the simulations, with loads controlled every 10ms. The disutility function for the deviation $d^c_j$ in controllable loads in each bus was $C_{dj}(d^c_j) = \frac{1}{2} \alpha_j ( d^{c}_j )^2$.  The selected values for cost coefficients were $\alpha_j = 1$ for load buses $1-5$ and $11-15$ and $\alpha_j = 2$ for the rest. Similarly, the cost functions for deviations $p^M_j$ in generation were $C_{j}(p^M_j) = \frac{1}{2} \kappa_j ( p^{M}_j )^2$, where $\kappa_j$ were selected as the inverse of the generators droop coefficients, as suggested in \eqref{contspec}.

Consider the static and first order dynamic schemes given~by $d^c_j = (C_{dj}')^{\hspace{-0.5pt}-1} (\omega_j - p^c_j)$ and
$\dot{d}^c_j = -d^c_j + (C_{dj}')^{\hspace{-0.5pt}-1}(\omega_j - p^c_j), j \in N$,
where $p^c_j$ has dynamics as described in \eqref{sys_power_command}. We refer to the resulting dynamics as Static and Dynamic OSLC respectively since in both cases, steady state conditions that solve the OSLC problem were {\color{black}used.}
{\color{black}As discussed in Section \ref{Discussion}, {\color{black}in} the presence of {\color{black}arbitrarily small} frequency damping,  both schemes satisfy Assumption \ref{assum2} and are thus included in our framework.}

The system was tested on three different cases. In case (i) 10 generators were employed to perform secondary frequency control by having frequency and power command as inputs. In case (ii) controllable loads were included on 20 load buses in addition to the 10 generators.
Controllable load dynamics in 10 buses were described by Static OSLC and in the rest by Dynamic OSLC. Finally, in case (iii), all controllable loads of case (ii) and {\color{black}15} generators where used for secondary frequency control. Note that the {\color{black}15} generators used for secondary frequency control had third, fourth and fifth order turbine governor dynamics.

The {\color{black}frequency at} bus 103 for the three tested cases {\color{black}is}  shown in Fig. \ref{Frequency}.
From {\color{black}this figure}, we observe that in all cases the frequency returns to its nominal value. However, the presence of controllable loads makes the frequency return much faster and with a smaller overshoot.

Furthermore, from Fig. \ref{Marginal}, it is observed
that the marginal costs at all controlled loads and generators that contribute to secondary frequency control, converge to the same value.  This illustrates the optimality in the power allocation 
among generators and loads, since equality in the marginal cost is necessary to solve \eqref{Problem_To_Min} {\color{black}when the power generated does not saturate to its maximum/minimum value}.}

\section{Conclusion}\label{Conclusion}
{\color{black}We have considered the problem of designing distributed schemes for secondary frequency control such that stability and optimality of the power allocation can be guaranteed. In particular, we have considered general classes of generation and demand control dynamics and have shown that a dissipativity condition in conjunction with appropriate decentralized conditions on their steady state {\color{black}behavior} can provide such stability and optimality guarantees.  We have also discussed that for linear systems the dissipativity condition can be easily verified by solving a corresponding LMI
{\color{black}and shown that the requirement to have knowledge of demand may be relaxed by incorporating an appropriate observer.}
Our results have been illustrated with simulations on the {\color{black}NPCC 140-bus} system.} Interesting potential extensions in the analysis include incorporating voltage dynamics, more advanced communication structures, as well as more advanced models for the loads where their switching {\color{black}behavior} is taken into account.
%

\vspace{-1mm}
\section*{{\color{black}Appendix A}} \label{Appendix}

{\color{black}
In this appendix we prove our main results, Theorems ~\ref{optthm} - \ref{convthm}, and also {\color{black}Lemmas \ref{eqbr_lemma}-\ref{eqbr_lemma_observer}}, Theorem \ref{suff_thm} and Proposition \ref{convthm_observer}.}

{\color{black}{\color{black}Throughout the proofs we will make use of the following equilibrium equations for the dynamics {\color{black}in \eqref{sys1}--\eqref{sys2},}
\begin{subequations} \label{eqbr}
\begin{equation}
0 = \omega^*_i - \omega^*_j, \; (i,j) \in E, \label{eqbr1}
\end{equation}
\begin{equation}
0 = - p_j^L + p_j^{M,*} - (d^{c,*}_j + d^{u,*}_j) - \sum_{k:j\rightarrow k} p^*_{jk} + \sum_{i:i\rightarrow j} p^*_{ij}, \; j\in G, \label{eqbr2}
\end{equation}
\begin{equation}
0 = - p_j^L - (d^{c,*}_j + d^{u,*}_j) - \sum_{k:j\rightarrow k} p^*_{jk} + \sum_{i:i\rightarrow j} p^*_{ij}, \; j\in L, \label{eqbr3}
\end{equation}
\begin{equation}
p^{M,*}_j = k_{p^M_j} (\zeta_j^*), \; j \in G, \label{eqbr5}
\end{equation}
\begin{equation}\label{eqbr6}
{\color{black}d^{c,*}_j = k_{d^c_j} (\zeta_j^*), \quad  
\zeta_j^{\ast} = [ -\omega_j^* \; p_j^{c,*}]^T, \quad j \in N. }
\end{equation}

\end{subequations}}

{\color{black}
\emph{Proof of Lemma~\ref{eqbr_lemma}:}
In order to show that $\omega^{*} = \vect{0}_{|N|}$, we sum equations \eqref{sys_power_command_b} at equilibrium for all $j \in N$, resulting in
$
\sum\limits_{j \in N} s^*_j = \sum\limits_{j \in N} p^L_j
$,
which shows that $\sum\limits_{j \in N} d^{u,*}_j = 0$ (by summing \eqref{eqbr2} and \eqref{eqbr3} over all $j \in G$ and $j \in L$ respectively). Then, Assumption \ref{assum2.4} implies that this equality holds only  if $\omega^{*} = \vect{0}_{|N|}$.
\hfill\IEEEQED
}

\emph{Proof of Theorem~\ref{optthm}:}
Due to Assumption~\ref{assum5}, $C_j'$ and $C_{dj}'$ are strictly increasing and hence invertible. Therefore all variables in \eqref{contspec} are
well-defined.
Also,  Assumption~\ref{assum5} guarantees that the OSLC {\color{black}optimization} problem~\eqref{Problem_To_Min} is convex and has a continuously differentiable cost function.
 Thus, a point $(\bar{p}^M, \bar{d}^c)$ is a global minimum for~\eqref{Problem_To_Min} if and only if it satisfies the KKT conditions~\cite{Boyd}
\begin{subequations} \label{kkt}
\begin{gather}
C_j'(\bar{p}^M_j) = \nu - \lambda_j^+ + \lambda_j^-, \; j \in G, \label{kkt1} \\
C_{dj}'(\bar{d}^c_j) = -\nu - \mu_j^+ + \mu_j^-, \; j \in N, \label{kkt2} \\
\sum\limits_{j\in  G} \bar{p}_j^M = \sum\limits_{j\in  N} (\bar{d}^c_j + p_j^L), \label{kkt4} \\
p^{M,min}_j \leq \bar{p}^M_j \leq p^{M,max}_j, \; j \in G, \label{kkt5} \\
d^{c,min}_j \leq \bar{d}^c_j \leq d^{c,max}_j, \; j \in N,  \label{kkt6} \\
\lambda_j^+ (\bar{p}^M_j - p^{M,max}_j) = 0, \; \lambda_j^- (\bar{p}^M_j - p^{M,min}_j) = 0, \; j \in G, \label{kkt7} \\
\mu^+ (\bar{d}^c_j - d^{c,max}_j) = 0, \; \mu^- (\bar{d}^c_j - d^{c,min}_j) = 0, \; j \in {\color{black}N}, \label{kkt8}
\end{gather}
\end{subequations}
for some constants $\nu \in \mathbb{R}$ and $\lambda_j^+, \lambda_j^-, \mu_j^+, \mu_j^- \ge 0$. It will be shown below that these conditions are satisfied by the equilibrium values $(\bar{p}^M, \bar{d}^c) = (p^{M,*}, d^{c,*})$ defined by {\color{black}equations~\eqref{eqbr5}, ~\eqref{eqbr6} and \eqref{contspec}.}

Since $C_j'$ and $C_{dj}'$ are strictly increasing, we can uniquely define
$\beta_j^{M,max} \!\! := \!\! C_j'(p^{M,max}_j)$, $\beta_j^{M,min} \!\! := \!\! C_j'(p^{M,min}_j)$, $\beta_j^{c,max} \!\! := -C_{dj}'(d^{c,max}_j)$, and $\beta_j^{c,min} \!\! := -C_{dj}'(d^{c,min}_j)$.  {\color{black}
 We let $\beta^*_0 = f(\zeta^{*}_j)$ {\color{black}where $f(.)$ is the function in the theorem statement. Note that} function $f$ is common
  {\color{black}at} every bus and is surjective {\color{black}hence} $\forall \beta^*_0 \in \mathbb{R}, \exists \zeta$ such that $f(\zeta) = \beta^*_0$.} {\color{black}Also note that the $\zeta_j$ are equal $\forall j$ at equilibrium, therefore $\beta^*_0$ is the same at each bus $j$.}
{\color{black}We now define in terms of these quantities} the nonnegative constants
\vspace{-1mm}
\begin{equation*}
\begin{aligned}
&\lambda_j^+ := ( \beta^*_0 - \beta^{M,max}_j) \, \mathds{1}_{(\beta^*_0 \ge \beta^{M,max}_j)}, \\
&\lambda_j^- := ( \beta^{M,min}_j - \beta^*_0 ) \, \mathds{1}_{(\beta^*_0 \le \beta^{M,min}_j)}, \\
&\mu_j^+ := ( \beta^{c,max}_j - \beta^*_0) \, \mathds{1}_{(\beta^*_0 \le \beta^{c,max}_j)}, \\
&\mu_j^- := ( \beta^*_0 - \beta^{c,min}_j) \, \mathds{1}_{(\beta^*_0 \ge \beta^{c,min}_j)}.
\end{aligned} \vspace{-1mm}
\end{equation*}
Then, since $(C_j')^{\hspace{-0.5pt}-1}(\beta_0^*) \ge p^{M,max}_j \Leftrightarrow \beta^*_0 \ge \beta^{M,max}_j$, $(C_j')^{\hspace{-0.5pt}-1}(\beta_0^*) \le p^{M,min}_j \Leftrightarrow \beta^*_0 \le \beta^{M,min}_j$, $(C_{dj}')^{\hspace{-0.5pt}-1}(-\beta_0^*) \ge d^{c,max}_j \Leftrightarrow \beta^*_0 \le \beta^{c,max}_j$, and $(C_{dj}')^{\hspace{-0.5pt}-1}(-\beta_0^*) \le d^{c,min}_j \Leftrightarrow \beta^*_0 \ge \beta^{c,min}_j$, it follows by~\eqref{eqbr5},~\eqref{eqbr6}, and~\eqref{contspec} that the complementary slackness conditions~\eqref{kkt7} and~\eqref{kkt8} are satisfied.

Now define $\nu = \beta_0^*$. Then $(C_j')^{\hspace{-0.5pt}-1} ( \nu - \lambda_j^+ + \lambda_j^- ) = (C_j')^{\hspace{-0.5pt}-1} \Big( [\beta_0^*]_{\beta^{M,min}_j}^{\beta^{M,max}_j} \Big) = [ (C_j')^{\hspace{-0.5pt}-1} (\beta_0^*)]_{p^{M,min}_j}^{p^{M,max}_j} = p^{M,*}_j$, by the above definitions and equations~\eqref{eqbr5} and~\eqref{contspec}. Thus, the optimality condition~\eqref{kkt1} holds. Analogously, $(C_{dj}')^{\hspace{-0.5pt}-1} ( -\nu - \mu^+ + \mu^- ) = (C_{dj}')^{\hspace{-0.5pt}-1} \Big( [-\beta_0^*]_{-\beta^{c,min}_j}^{-\beta^{c,max}_j} \Big) = [ (C_{dj}')^{\hspace{-0.5pt}-1} (-\beta_0^*)]_{d^{c,min}_j}^{d^{c,max}_j} = d^{c,*}_j$, by~\eqref{eqbr6} and~\eqref{contspec}, satisfying~\eqref{kkt2}.

\balance
Summing {\color{black} equations ~\eqref{eqbr2} and~\eqref{eqbr3} over all $j \in G$ and $j \in L$ respectively} and using the fact {\color{black}that $\sum_{j \in N} d^{u,*}_j = 0$ as shown in the proof of Lemma \ref{eqbr_lemma}} shows that~\eqref{kkt4} holds. Finally, the saturation constraints  in~\eqref{contspec} verify~\eqref{kkt5}~and~\eqref{kkt6}.

{\color{black}Hence}, the values $(\bar{p}^M, \bar{d}^c) = (p^{M,*}, d^{c,*})$ satisfy the KKT conditions~\eqref{kkt}. Therefore, the equilibrium values $p^{M,*}$ and $d^{c,*}$ define a global minimum for~\eqref{Problem_To_Min}. \hfill\IEEEQED

\emph{Proof of Theorem~\ref{convthm}:}
We will use the dynamics in \eqref{sys1}--\eqref{sys_power_command} and the conditions of Assumption~\ref{assum2} to define a Lyapunov function for the system~\eqref{sys1}--\eqref{sys_power_command}.

Firstly, let $V_F (\omega^G) = \frac{1}{2}\sum_{j \in G} M_j (\omega_j - \omega^*_j)^2$. {\color{black} The time-derivative of $V_F$} along the trajectories of~\eqref{sys1}--\eqref{sys2} {\color{black}is} given by
\begin{align*}
\hspace{-0mm}\dot{V}_F &=\hspace{-1mm} \sum_{j \in N}\hspace{-0.5mm}(\omega_j - \omega^*_j) \Bigg(\hspace{-0.5mm}-p^L_j + s_j - d^u_j -\hspace{-1mm} \sum_{k:j\rightarrow k} p_{jk} + \hspace{-1mm}\sum_{i:i\rightarrow j} p_{ij}\hspace{-0.5mm}\Bigg),
\end{align*}
by substituting~\eqref{sys1b} for $\dot{\omega}_j$ for $j \in G$ and adding extra terms for $j \in L$, which are equal to zero by~\eqref{sys1c}. Subtracting the product of $(\omega_j - \omega^*_j)$ with each term in~\eqref{eqbr2} and~\eqref{eqbr3}, this becomes
\begin{align}
\dot{V}_F =& \sum_{j \in N} \Bigg ({\color{black}\hspace{-0.5mm}(\omega_j} - \omega^*_j) (s_j - s^{*}_j)
\hspace{-0.5mm} + \hspace{-0.2em} (\omega_j - \omega^*_j) (-d^u_j - (-d^{u,*}_j))\hspace{-0.5mm}\Bigg)  \nonumber \\
& {\color{black}+
\sum_{(i,j) \in E}} (p_{ij} - p^*_{ij}) (\omega_j - \omega_i), \label{VFdiff}
\end{align}
using the equilibrium condition~\eqref{eqbr1} for the final term.

Furthermore, let $V_C(p^c) = \frac{1}{2}\sum_{j \in N} \gamma_j(p^c_j - p^{c,*}_j)^2$. {\color{black}Using \eqref{sys_power_command_b} the} time derivative of $V_C$ can be written as
\begin{multline}
\dot{V}_C = \sum_{j \in N}  (p^c_j - p^{c,*}_j) \Big((-s_j + s^{*}_j) \\
- \sum_{k:j\rightarrow k} (\psi_{jk} - \psi^*_{jk}) + \sum_{i:i\rightarrow j} (\psi_{ij} - \psi^*_{ij})\Big).
\end{multline}

Additionally, define $V_P(\eta) = \sum_{(i,j) \in E} B_{ij} \int_{\eta^*_{ij}}^{\eta_{ij}} ( \sin \theta - {\color{black}\sin \eta^*_{ij}} ) \, d\theta$. Using~\eqref{sys1a} and~\eqref{sys1d}, the time-derivative {\color{black}is given~by}
\begin{align}
\dot{V}_P &= \sum_{(i,j) \in E} B_{ij} (\sin \eta_{ij} - \sin \eta^*_{ij}) (\omega_i - \omega_j) \nonumber \\
&= \sum_{(i,j) \in E} (p_{ij} - p^*_{ij}) (\omega_i - \omega_j). \label{VPdiff}
\end{align}

Finally, consider {\color{black}$V_{\psi}(\psi) = \frac{1}{2}  \sum_{(i,j) \in {\color{black}\tilde E}}\gamma_{ij}(\psi_{ij} - \psi^*_{ij})^2$} with time derivative given by \eqref{sys_power_command_a} as

\begin{equation}\label{Vphi_diff}
\dot{V}_{\psi} = \sum_{(i,j) \in {\color{black}\tilde E}} (\psi_{ij} - \psi^*_{ij})((p^c_i - p^{c,*}_i) - (p^c_j - p^{c,*}_j)).
\end{equation}

Furthermore, from the dissipativity {\color{black}condition in Assumption~\ref{assum2}
the following holds:}
{\color{black}
There exist open neighborhoods $U_j$ of $\omega^*_j$ and {\color{black}$U^c_j$} of $p^{c,*}_j$ for each $j \in N$, open {\color{black}neighborhoods} $X^G_j$ of $(x^{M,j,*},x^{c,j,*}, x^{u,j,*})$ and $X^L_j$ of $(x^{c,j,*}, x^{u,j,*})$  for each $j \in G$ and $j \in L$ respectively,
 and continuously differentiable, positive semidefinite {\color{black}functions}
 {\color{black}$V^D_j(x^{M,j},x^{c,j}, x^{u,j}), j \in G$  and $V^D_j(x^{c,j}, x^{u,j}), j \in L$,} satisfying \eqref{dissipativity_condition} with {\color{black}supply rate given by \eqref{supply_rate}, i.e.,}
\begin{multline}\label{VG_derivative}
\dot{V}^D_j \le
[(s_j - s^{*}_j) \;\;\;  (-d^u_j - (-d^{u,*}_j))] \begin{bmatrix}
1 & 1\\
1 & 0
\end{bmatrix} (\zeta_j - \zeta^*_j) \\
{\color{black}- {\color{black}\phi_j}(\zeta_j - \zeta^*_j),}
\hspace{3mm}
j \in N,
\end{multline}
for all $\omega_j \in U_j$, $p^{c}_j$ in $U^c_j$ for $j \in N$ and all $(x^{M,j},x^{c,j}, x^{u,j}) \in X^G_j$ and $(x^{c,j}, x^{u,j}) \in X^L_j$ for $j \in G$ and $j \in L$ respectively.}

{\color{black}Based on the above, we define the function
{\color{black}
\begin{equation}\label{Lyapunov_function}
V(\eta, \psi, \omega^G, x^M, x^c, x^u, p^c) = V_F + V_P
 + \sum_{j \in N} V^D_j
 + V_C + V_{\psi}
\end{equation}
which} we aim to use in Lasalle's theorem.
Using \eqref{VFdiff}  - \eqref{Vphi_diff}, the time derivative of V is {\color{black}given by
\begin{align}\label{dot_V}
\dot{V} =& \sum_{j \in N} \bigl[(\omega_j - \omega^*_j) (s_j - s^{*}_j)
 + \dot{V}^D_j
+ (p^c_j - p^{c,*}_j) (-s_j + s^{*}_j) \nonumber\\
&\hspace{2.5cm} + (\omega_j - \omega^*_j) (-d^u_j - (-d^{u,*}_j)\bigr].
\end{align}}
Using~\eqref{VG_derivative}
it therefore {\color{black}holds that
\begin{align}
\dot{V} &\le \sum_{j \in N} \Big(
- {\color{black}\phi_j}(\zeta_j -{\color{black} \zeta^*_j)}
\Big) \le0
\label{vdotineq}
\end{align}
whenever $\omega_j \in U_j$, $p^c_j \in {\color{black}U^c_j}$ for $j \in N$, $(x^{M,j},x^{c,j}, x^{u,j}) \in X^G_j$ for $j \in G$, and $(x^{c,j},x^{u,j}) \in X^L_j$ for $j \in L$.} 

Clearly $V_F$ has a strict global minimum at $\omega^{G,*}$ and {\color{black} $V^D_j$ has strict local minima at $(x^{M,j,*}, x^{c,j,*},x^{u,j,*})$ and $(x^{c,j,*},x^{u,j,*})$ for $j \in G$ and $j \in L$ respectively by Assumption~\ref{assum2} and Definition \ref{Dissipativity_Definition}.} Furthermore, $V_C$ and $V_{\psi}$ have strict global minima at $p^{c,*}$ and $\psi^*$ respectively. Furthermore, Assumption~\ref{assum1} guarantees the existence of some {\color{black}neighborhood} of each $\eta^*_{ij}$ {\color{black}in which} $V_P$ is increasing. Since the integrand is zero at the lower limit of the integration, $\eta^*_{ij}$, this immediately implies that $V_P$ has a strict local minimum at $\eta^*$. Thus, $V$ has a strict local minimum at the point $Q^* := (\eta^*,\psi^*, \omega^{G,*}, x^{M,*}, x^{c,*}, x^{u,*}$, $p^{c,*})$.
From Assumption~\ref{assum4}, we know that, provided {\color{black}$(\eta, \omega^G, x^c, x^u, p^c) \in T$,} $\omega^L$ can be uniquely determined from these quantities. Therefore, the states of the differential equation system~\eqref{sys1}--\eqref{sys_power_command} with {\color{black}$(\eta, \omega^G, x^c, x^u, p^c)$  within the region $T$} can be expressed as $(\eta, \psi, \omega^G, x^M, x^c, x^u, p^c)$. We now choose a {\color{black}neighborhood} in the coordinates $(\eta, \psi, \omega^G, x^M, x^c, x^u, p^c)$ about $Q^*$ on which the following hold:
\begin{enumerate}
\item $Q^*$ is a strict minimum of $V$,
\item {\color{black}$(\eta, \omega^G, x^c, x^u, p^c) \in T$,}
\item {\color{black}$\omega_j \in U_j$, $p^c_j \in {\color{black}U^c_j}$ for $j \in N$, and  $(x^{M,j},x^{c,j},x^{u,j}) \in X^G_j$, $(x^{c,j},x^{u,j}) \in X^L_j$ for $j \in G$, $j \in L$ respectively\footnote{This is possible {\color{black}because} {\color{black}$\omega_j \in U_j$} for all $j \in L$ corresponds, by Assumption~\ref{assum4} and the continuity of the {\color{black}equations in~\eqref{sys1}--\eqref{sys_power_command}}, to requiring the states {\color{black}$(\eta, \omega^G, x^M, x^c, x^u,p^c)$} to lie in some open {\color{black}neighborhood} about $Q^*$.},}
\item $x^{M,j}$, $x^{c,j}$, and $x^{u,j}$ all lie within their respective {\color{black}neighborhoods} $X_0$ as defined in {\color{black}Section~\ref{Network_Model}}.
\end{enumerate}
Recalling now~\eqref{vdotineq}, it is easy to see that within this {\color{black}neighborhood,} $V$ is a nonincreasing function of all the system states and has a strict local minimum at $Q^*$. Consequently, the connected component of the level set $\{(\eta, \psi, \omega^G, x^M, x^c, x^u, p^c) \colon V \le \epsilon\}$ containing $Q^*$ is guaranteed to be both compact and positively invariant with respect to the system~\eqref{sys1}--\eqref{sys_power_command} {\color{black}for sufficiently} small $\epsilon > 0$. Therefore, there exists a compact positively invariant set $\Xi$ for~\eqref{sys1}--\eqref{sys_power_command} containing~$Q^*$.

Lasalle's Invariance Principle can now be applied with the function $V$ on the compact positively invariant set~$\Xi$. This guarantees that all solutions of~\eqref{sys1}--\eqref{sys_power_command} with initial conditions $(\eta(0), \psi(0), \omega^G(0), x^M(0), x^c(0) ,x^u(0), p^c(0)) \in \Xi$ converge to the largest invariant set within $\Xi \, \cap \, \{(\eta, \psi, \omega^G, x^M, x^c, x^u, p^c) \colon \dot{V} = 0\}$. We now consider this invariant set. If $\dot{V} = 0$ holds at a point within $\Xi$, then~\eqref{vdotineq} holds with equality, hence  we must have $\omega = \omega^*$ and $p^c_j = p^{c,*}_j$
{\color{black}
at all buses $j$ where Assumption \ref{assum2}(a) holds.  The fact that $\omega$ is constant guarantees from \eqref{sys1a}, \eqref{sys1d} that $\eta$ and $p$ are also constant. This is sufficient to deduce from {\color{black} \eqref{sys1b}--\eqref{sys1c}} that $s$ is also constant.
If instead Assumption \ref{assum2}(b) holds at a bus $j$ we have {\color{black} that $\omega = \omega^*$ when $\dot V =0$}. Furthermore, we have the additional property that {if  $\omega_j$ and $s_j$ are constant then $p^c_j$ cannot be a sinusoid. This latter property guarantees that $p^c_j$ is also  constant by
noting that the dynamics for the power command \eqref{sys_power_command} with constant $s_j$,  allow $p^c_j$ to be either a constant or a sinusoid within a compact invariant set.}
Hence, we have $\omega = \omega^*$ and $p^c = p^{c,*}$ in the invariant set considered.}

Furthermore, note that $\omega = \omega^*$, $p^c = p^{c,*}$ within the invariant set
implies by the definitions in Section~\ref{sec: Preliminaries} that $(x^M, x^c, x^u)$ converge to the point $(x^{M,*}, x^{c,*}, x^{u,*})$, at which  {\color{black}$V^D_j$} take strict local minima from Assumption~\ref{assum2}.  Thus,
from {\color{black}~\eqref{VG_derivative} and \eqref{dot_V}  it follows that the values of $V^D_j$ must decrease along all nontrivial trajectories within the invariant set, contradicting $\dot{V}^D_j = 0$.}
 The fact that $(p^c,s)= (p^{c,*},s^*$) is sufficient to show that $\psi$ equals some constant $\psi^*$.
 Using the same argument, {\color{black} it can be shown that within the invariant set,
 {\color{black}the fact that $\zeta = \zeta^*$ implies that} $(x^M,x^c,x^u,p^M,d^c,d^u)$ converges to $(x^{M,*},x^{c,*},x^{u,*},p^{M,*},d^{c,*},d^{u,*})$.} Therefore, we conclude by Lasalle's Invariance Principle that all solutions of~\eqref{sys1}--\eqref{sys_power_command} with initial conditions $(\eta(0), \psi(0), \omega^G(0), x^M(0), x^c(0) ,x^u(0), p^c (0)) \in \Xi$ converge to the set of equilibrium points as defined in Definition~\ref{eqbrdef}. Finally, choosing for $S$ any open {\color{black}neighborhood} of $Q^*$ within $\Xi$ completes the proof for convergence. {\color{black}{\color{black}From} Lemma \ref{eqbr_lemma} it can then be deduced that $\omega^{*} = \vect{0}_{|N|}$. {\color{black}Furthermore, noting that all conditions of Theorem \ref{optthm} hold shows the convergence to an optimal solution of the OSLC problem \eqref{Problem_To_Min}. }}
\hfill\IEEEQED

{\color{black}
\begin{remark}
It should be noted that for given $p^{c,*}$ and $\omega^{*}$ all $(\eta^*, x^{M,*},x^{c,*},x^{u,*})$ are unique. The uniqueness of $\eta^*$ can be seen by noting that $\eta_{ij} = \theta_i - \theta_j, (i,j) \in E$, which requires $\eta$ to lie in a space where a corresponding vector $\theta$ exists.
Furthermore, the value of $p^{c,*}$ becomes unique when \eqref{contspec} holds. This follows from summing \eqref{eqbr2}--\eqref{eqbr3} over all buses and noting that the strict convexity of the cost functions {\color{black}and the monotonicity of $f$ in \eqref{contspec}} makes the static input output maps from $p^{c,*}$ to $s^*$ monotonically increasing.
The values of $\psi^*$ are non-unique for general network topologies.
\end{remark}
}

{\color{black}{\color{black}
\emph{Proof of Theorem~\ref{suff_thm}:}
The proof follows from the fact that the function $V^B_j$ defined as
\vspace{-.2cm}
\begin{multline}\label{Storage_function}
V^B_j = \frac{1}{2}  M_j (\omega_j - \omega^*_j)^2 +
 \frac{1}{2} \gamma_j(p^c_j - p^{c,*}_j)^2 +
 V^D_j, 
\end{multline}
where $V^D_j$ is as in \eqref{VG_derivative}  with $\phi_j=0$, is a storage function for the system $B_j$.
In particular,
using arguments similar to those in the proof of Theorem \ref{convthm}, it can be shown that
\begin{multline}\label{Storage_derivative}
\dot{V}^B_j \leq (p^c_j - p^{c,*}_j)\biggl( \sum_{i:i\rightarrow j} (\psi_{ij}-\psi_{ij}^\ast) -\sum_{k:j\rightarrow k} (\psi_{jk
}-\psi_{jk}^\ast)\biggr) \\
+ (- \omega_j - (-\omega^*_j))\biggl(\sum_{k:j\rightarrow k} (p_{jk}-p_{jk}^\ast) - \sum_{i:i\rightarrow j} (p_{ij}-p_{ij}^\ast)\biggr)
\end{multline}
{\color{black}and therefore that} system $B_j$ is passive.
\hfill\IEEEQED}
}

%

{\color{black}
\emph{Proof of Lemma~\ref{eqbr_lemma_observer}:}
Using \eqref{sys_observer_d} at equilibrium, it can be deduced that $\chi^*_j =  s^*_j  - d^{u,*}_j - \sum_{k:j\rightarrow k} p^*_{jk} + \sum_{i:i\rightarrow j} p^*_{ij}$. Hence, it follows by summing \eqref{sys_observer_b} at equilibrium over all buses  that  $\sum\limits_{j \in N} s^*_j = \sum\limits_{j \in N} \chi^*_j = \sum\limits_{j \in N} s^*_j - d^{u,*}_j $, which results to $\sum \limits_{j \in N} d^{u,*}_j  = 0$. Hence, from Assumption \ref{assum2.4}, it follows that $\omega^* = \vect{0}_{|N|}$.
\hfill\IEEEQED

\emph{Proof of Proposition~\ref{convthm_observer}:}
We shall make use of the Lyapunov function in \eqref{Lyapunov_function}  to construct a new Lyapunov function for the system \eqref{sys1} -- \eqref{ssys}, \eqref{sys_observer}.

First, consider the function
\begin{multline}\label{Vb}
V_b(b, \chi, \omega) =  \frac{1}{2}\sum_{j \in G} \big( M_j ((b_j - b^*_j) - (\omega_j - \omega^*_j))^2 \\
 + \tau_{\chi,j} (\chi_j - \chi^*_j)^2 \big),
\end{multline}
and note that its time-derivative along the trajectories of \eqref{sys_observer} is given by
\begin{equation}\label{dot_Vb}
\dot{V}_b = \sum_{j \in N} \biggr( -(\chi_j - \chi^*_j) [(p^c_j - p^{c,*}_j) + (\chi_j - \chi^*_j)] \biggr),
\end{equation}
noting that for $j \in L$ it holds that $\chi = \chi^*$, and hence the added terms in \eqref{dot_Vb} are equal to zero.

%

Furthermore, the time-derivative of $V_C(p^c) = \frac{1}{2}\sum_{j \in N} \gamma_j(p^c_j - p^{c,*}_j)^2$ under \eqref{sys_observer_b} is given by
\begin{multline}\label{dot_VC}
\dot{V}_C = \sum_{j \in N}  (p^c_j - p^{c,*}_j) \Big((-s_j + s^{*}_j) +  (\chi_j - \chi^*_j) \\
- \sum_{k:j\rightarrow k} (\psi_{jk} - \psi^*_{jk})
+ \sum_{i:i\rightarrow j} (\psi_{ij} - \psi^*_{ij}) \Big).
\end{multline}
 Now consider the function $V$ in \eqref{Lyapunov_function} and note that its derivative is as in \eqref{dot_V} with an extra term given by $\sum_{j \in N}  (p^c_j - p^{c,*}_j)(\chi_j - \chi^*_j)$. Then consider the function
\begin{equation}
V_O(\eta, \psi, \omega^G, x^M, x^c, x^u, p^c, b, \chi) = V + V_b
\end{equation}
which can be shown to have a time derivative given by
\begin{align}
\dot{V}_O &\le  \sum_{j \in N}
\Big(-\phi_j(\zeta_j - \zeta^*_j) - (\chi_j - \chi^*_j)^2 \Big)
 \le0,
\label{vdotineq_observer}
\end{align}
by similar arguments as in the proof of Theorem \ref{convthm}.

Now, in analogy to the proof of Theorem \ref{convthm}, it can be shown that an invariant compact set $\Xi_O$ exists such that  $\{(\eta, \psi, \omega^G, x^M, x^c, x^u, p^c, b, \chi) \colon V_O \le \epsilon\}$. Then, Lasalle's theorem can be invoked to show that all solutions of \eqref{sys1} -- \eqref{ssys}, \eqref{sys_observer} with initial conditions within $\Xi_O$ will converge to the largest invariant set within $\Xi_O \, \cap \, \{(\eta, \psi, \omega^G, x^M, x^c, x^u, p^c,  b, \chi) \colon \dot{V} = 0\}$. Within this invariant set, it holds that
 $(\omega, \chi) = (\omega^*, \chi^*)$. Applying the same arguments as in the proof of Theorem \ref{convthm} shows that $(\eta, \psi, x^M,x^c,x^u,p^M,d^c,d^u,p^c)$ converges to $(\eta^*, \psi^*, x^{M,*},x^{c,*},x^{u,*},p^{M,*},d^{c,*},d^{u,*}, p^{c,*})$ which implies the convergence of $b$ to $b^*$ from the dynamics in \eqref{sys_observer_d}.
 {\color{black}The optimality result follows directly from the proof of Theorem \ref{optthm} since none of its arguments are affected from the dynamics in \eqref{sys_observer}.}
\hfill\IEEEQED
}

{\color{black}
\section*{Appendix B}\label{Appendix_B}

In this appendix we show that Assumption \ref{assum2} is a necessary {\color{black}and sufficient} condition for the passivity of bus systems $B_j$, described in Section \ref{System_Representation}, when their dynamics are affine nonlinear, i.e. are characterized by the following state space representation:
\begin{align}\label{affine_systems}
\dot{x} = f(x) + g(x)u,  \nonumber \\
 y = h(x).
\end{align}

{\color{black}
For the proof, we shall make use of Lemma \ref{Sufficiency_lemma} below.
Within it, we shall consider the negative feedback interconnection of}
\begin{equation}\label{Sigma_feed}
\Sigma_1:
\begin{cases}
\dot x_1=u_1\\
y_1=h_1(x_1)
\end{cases},
\Sigma_2:
\begin{cases}
\dot x_2=f_2(x_2)+g_2(x_2)u_2\\
y_2=h_2(x_2)+k_2(u_2)
\end{cases}
\end{equation}
such that
$u_2 = y_1$ and $u_1=r-y_2$,
where $r(t) \in \mathbb{R}^n$ is some reference input applied to the closed-loop system,  $x_1(t) \in \mathbb{R}^n, x_2(t) \in \mathbb{R}^{n_2}$ and $y_1(t), y_2(t) \in \mathbb{R}^n$  are the states and outputs of $\Sigma_1$ and $\Sigma_2$ respectively and $h_1,k_2 : \mathbb{R}^n \rightarrow \mathbb{R}^n$, and $f_2,g_2, h_2: \mathbb{R}^{n_2} \rightarrow \mathbb{R}^n$ are functions describing the {\color{black}dynamics} of $\Sigma_1$ {\color{black}and $\Sigma_2$.} 
}
The closed-loop system, denoted by $\Sigma$, writes as
\begin{align}\label{Sigma}
\Sigma: \begin{cases}\dot x_1&=-h_2(x_2)-k_2(h_1(x_1))+r,\\
\dot x_2&=f_2(x_2)+g_2(x_2)h_1(x_1),\\
y_1&=h_1(x_1).
\end{cases}
\end{align}
{\color{black}Without loss of generality we also assume in the {\color{black}lemma} below that $h_1(0),k_2(0),f_2(0),h_2(0)$ are equal to zero and the passivity properties stated are considered about this {\color{black}equilibrium point.}}


\begin{lemma}\label{Sufficiency_lemma}
Consider the negative feedback interconnection $\Sigma$ described by \eqref{Sigma} of two subsystems $\Sigma_1$ and $\Sigma_2$  described by \eqref{Sigma_feed}.
Assume that $\Sigma$ is passive from $r$ to $y_1$, and $\Sigma_1$ is passive from $u_1$ to $y_1$.
Then, $\Sigma_2$ is passive from $u_2$ to $y_2$.
\end{lemma}
\emph{Proof of Lemma~\ref{Sufficiency_lemma}:}
From the passivity of $\Sigma$ and \cite[Corollary 4.1.5]{L2_gain_passivity} there exists a positive definite continuously differentiable storage function $V(x_1,x_2)$, defined with respect to an equilibrium, such that
\begin{align}\label{main}
-\frac{\partial{V}}{\partial{x_1}}(x_1,x_2)h_2(x_2)
-\frac{\partial{V}}{\partial{x_1}}(x_1,x_2)k_2(h_1(x_1)) \nonumber\\
+\frac{\partial{V}}{\partial{x_2}}(x_1,x_2)f_2(x_2)
+\frac{\partial{V}}{\partial{x_2}}(x_1,x_2)g_2(x_2)h_1(x_1)
\leq 0
\end{align}
and
\[
\frac{\partial{V}}{\partial{x_1}}(x_1, x_2)=h_1^T(x_1).
\]
Similarly, from the passivity of $\Sigma_1$ and \cite[Corollary 4.1.5]{L2_gain_passivity} there exists a positive definite continuously differentiable storage function $V_1$ such that
\begin{align}
\label{sig1-x1}
\frac{\partial{V_1}}{\partial{x_1}}(x_1)=h_1^T(x_1).
\end{align}
Hence,
\begin{equation}\label{e:imp}
\frac{\partial{V}}{\partial{x_1}}(x_1,x_2)=\frac{\partial{V_1}}{\partial{x_1}}(x_1).
\end{equation}
Substituting this back to \eqref{main} yields
\begin{align}\label{main2}
\nonumber
\frac{\partial{V}}{\partial{x_2}}&(x_1,x_2)f_2(x_2)
+
\frac{\partial{V}}{\partial{x_2}}(x_1,x_2)g_2(x_2)h_1(x_1)\\
&\leq
\frac{\partial{V_1}}{\partial{x_1}}(x_1)(h_2(x_2)+k_2(h_1(x_1))  \nonumber\\
&= h_1^T(x_1)(h_2(x_2)+k_2(h_1(x_1)),
\end{align}
where the last inequality follows from \eqref{sig1-x1}.
Now let $V$ be written as
\begin{equation}\label{V_summation}
V(x_1, x_2)=V_1(x_1)+V_2(x_2)
\end{equation}
for some continuously differentiable $V_2$. The fact that $V_2$ is only a function of $x_2$ descends from \eqref{e:imp}.
Also note that $V_2(x_2)$ is positive definite.
By substituting the above into \eqref{main2}, we conclude that
\[
\frac{\partial{V_2}}{\partial{x_2}}(x_2)f_2(x_2)
+
\frac{\partial{V_2}}{\partial{x_2}}(x_2)g_2(x_2)u_2\leq
u^T_2y_2\]
which implies the passivity
of $\Sigma_2$.
\hfill \IEEEQED

The following lemma shows that Assumption \ref{assum2} is a necessary {\color{black}and sufficient} condition for the passivity of {\color{black}generation} bus system $B_j$. Note that the extension to load buses is trivial and thus omitted.

\begin{lemma}\label{suff_on_Bj}
Consider the system described by \eqref{sys1} - \eqref{sys_power_command} and its representation by systems $I$ and $B$, defined in section \ref{System_Representation} and let the dynamics for $B_j$ be described {\color{black}by \eqref{Sigma}.}
Then, the dissipativity condition in Assumption \ref{assum2}  with $\phi=0$ {\color{black}is necessary and sufficient} for the passivity of subsystems $B_j, j\in G$ about the equilibrium point considered in Assumption \ref{assum2}.
\end{lemma}

\emph{Proof of Lemma~\ref{suff_on_Bj}:}
{\color{black}The proof for the necessity of the condition} follows from Lemma \ref{Sufficiency_lemma} when the following substitutions are {\color{black}made
\begin{align*}
r =
\begin{bmatrix}
\sum_{k:j\rightarrow k} (p_{jk}-p^*_{jk}) - \sum_{i:i\rightarrow j} (p_{ij}-p^*_{ij})   \\
    \sum_{i:i\rightarrow j} (\psi_{ij} -\psi^*_{ij})  - \sum_{k:j\rightarrow k} (\psi_{jk}-\psi^*_{jk})
\end{bmatrix}, \nonumber \\
y_1 = \begin{bmatrix}
-(\omega_j - \omega^*_j) \\
p^c_j - p^{c,*}_j
\end{bmatrix}, \;
y_2 = \begin{bmatrix}
(s_j-s_j^\ast) - (d^u_j-d^{u,\ast}_j)\\
(s_j-s_j^\ast)
\end{bmatrix}.
\end{align*}
{\color{black}The sufficiency proof follows directly from Theorem \ref{suff_thm}.}}
\hfill \IEEEQED
}
\balance

\end{document}